\documentclass{amsart}

\usepackage{m-pictex, epsfig, amsthm, amssymb, amsmath, amsxtra, latexsym}
\usepackage[all]{xy}
\usepackage{enumerate}

\newcommand{\mZ}{\mathbb Z}
\newcommand{\tK}{\tilde K_0}

 \DeclareMathOperator{\z}{{\mathbb
Z}} 
\DeclareMathOperator{\g}{\Gamma}

\DeclareMathOperator{\kz}{\mathbb K \mathbb Z^{-\infty}}
\DeclareMathOperator{\vc}{\mathcal V\mathcal C}

\DeclareMathOperator{\fin}{\mathcal F\mathcal I\mathcal N}
\DeclareMathOperator{\all}{\mathcal A\mathcal L\mathcal L}

\theoremstyle{plain}
\newtheorem{theorem}{Theorem}[section]
\newtheorem{lemma}[theorem]{Lemma}

\newtheorem*{main}{Main Theorem}
\newtheorem*{Thm}{Corollary}

\theoremstyle{definition}

\theoremstyle{remark}

\begin{document}

\title[Relating the Farrell Nil-groups to the Waldhausen Nil-groups.]
      {Relating the Farrell Nil-groups to the Waldhausen Nil-groups.}
\author{Jean-Fran\c{c}ois\ Lafont}
\address{Department of Mathematics\\
         Ohio State University\\
         Columbus, OH  43210}
\email[Jean-Fran\c{c}ois\ Lafont]{jlafont@math.ohio-state.edu}
\author{Ivonne J.\ Ortiz}
\address{Department of Mathematics and Statistics\\
         Miami University\\
         Oxford, OH 45056}
\email[Ivonne J.\ Ortiz]{ortizi@muohio.edu}

\begin{abstract}
Every virtually cyclic group $\g$ that surjects onto the infinite
dihedral group $D_\infty$ contains an index two subgroup $\Pi$ of
the form $H\rtimes _\alpha \mathbb Z$.  We show that the Waldhausen
Nil-group of $\g$ vanishes if and only if the Farrell Nil-group of
$\Pi$ vanishes.
\end{abstract}

\maketitle

\section{Statement of results.}

The Bass Nil-groups, Farrell Nil-groups, and Waldhausen Nil-groups
appear respectively as pieces in the computation of the algebraic
$K$-theory of direct products, semi-direct products, and
amalgamations.  While the Bass Nil-groups have been extensively
studied, much less is known for both the Farrell Nil-groups and the
Waldhausen Nil-groups.  For the purposes of computing the
algebraic $K$-theory of infinite groups, the Nil-groups of virtually
cyclic groups yield obstructions to certain assembly maps being
isomorphisms. In particular, the vanishing/non-vanishing of
Nil-groups is of crucial importance for computational aspects of
algebraic $K$-theory. In this short note we prove the
following result:

\begin{main}
Let $\g$ be a virtually cyclic group that surjects onto the infinite dihedral group
$D_\infty$, and $\g=G_1 \ast_{H} G_2$ be the corresponding
splitting of groups (with $H$ of index two in both $G_1$ and $G_2$).
Let $\Pi= H \rtimes_{\alpha} \mathbb Z\leq \g$ be the canonical
subgroup of $\g$ of index two, obtained by taking the pre-image of
the canonical index two $\mZ$ subgroup of $D_\infty$. Then for $i=0,1$, the
following two statements are equivalent:
\begin{itemize}
\item[(A)]  The Waldhausen Nil-group $NK_i(\mathbb ZH;\mathbb Z[G_1 - H],
\mathbb Z[G_2 - H])$ for the group $\g =G_1 \ast_{H} G_2$ vanishes.
\item[(B)]  The Farrell Nil-group $NK_i(\mathbb ZH, \alpha)$ for the group
$\Pi = H\rtimes_{\alpha} \mathbb Z$ vanishes.
\end{itemize}
\end{main}

The proof of our Main Theorem will be completed in Section 2, with
some concluding remarks in Section 3.  




Next, let us recall that the
Farrell-Jones Isomorphism Conjecture for a finitely generated group
$\g$ states that the assembly map:
$$H_n^{\g} (E_{\vc}(\g); \kz)\longrightarrow H_n^{\g}(E_{\all}(\g);\kz) =K_{n}(\mathbb Z\g)$$
is an isomorphism.  The term on the left is the generalized equivariant homology
theory of the space $E_{\vc}(\g)$ with coefficients in the integral
$K$-theory spectrum, where the space $E_{\vc}(\g)$ is a classifying
space for $\g$-actions with isotropy in the family $\vc$ of
virtually cyclic subgroups.  The term on the right gives the
algebraic $K$-theory of the integral group ring of $\g$.

Explicit models for the classifying space $E_{\vc}(\g)$ are known
for few classes of groups: virtually cyclic groups (take
$E_{\vc}(\g)$ to be a point with trivial action), crystallographic
groups (by work of Alves and Ontaneda \cite{AO06}), hyperbolic
groups (by work of Juan-Pineda and Leary \cite{JL}, and L\"uck
\cite{Lu05}), and in the case of relatively hyperbolic groups (due
to the authors \cite{LO}).  In contrast, explicit models are known
for $E_{\fin}(\g)$ for many classes of groups (see \cite{Lu05} for a
thorough survey), where $E_{\fin}(\g)$ is a classifying space for
$\g$-actions with isotropy in the family of finite subgroups.  Furthermore,
the algebraic K-theory of finite groups is much better understood than
the algebraic K-theory of virtually cyclic groups.  As
such, it is interesting to know whether one can further reduce the
computation of the algebraic $K$-theory of the integral group ring
of $\g$ to the computation of $H_n^{\g} (E_{\fin}(\g); \kz)$, i.e.,
whether the natural map:
$$H_n^{\g} (E_{\fin}(\g); \kz)\rightarrow H_n^{\g} (E_{\vc}(\g); \kz)$$
is an isomorphism.  There is a well known criterion for this, namely
the map will be an isomorphism for all $n \leq q$ if and only if,
for every infinite virtually cyclic subgroup $Q$ of $\g$, the map:
$$H_n^{Q} (E_{\fin}(Q); \kz)\rightarrow K_n(\z Q)$$
is an isomorphism for all $n \leq q$ (see \cite[Theorem 2.3]{DL98}
for $q<\infty$, and \cite[Theorem A.10]{FJ93} for $q=\infty$).  Since
the cokernel of these maps are precisely the Nil-groups, our main
theorem now gives us the following

\begin{Thm}
Let $\g$ be a finitely generated group, and $q=0$ or $1$.  Then the following
two statements are equivalent:
\begin{itemize}
\item  the relative assembly map
$$H_n^{\g} (E_{\fin}(\g); \kz)\rightarrow H_n^{\g}(E_{\vc}(\g);\kz)$$
is an isomorphism for $n\leq q$,
\item for every subgroup $H\rtimes _\alpha \mZ \leq \g$ with $H$ finite, and for
every  $n\leq q$, the
Farrell Nil-group $NK_n(\mathbb ZH, \alpha)$ for the group
$H\rtimes_\alpha \mZ$ vanishes.
\end{itemize}
If in addition, we know that the
Farrell-Jones Isomorphism Conjecture holds for the group $\g$, then
we have that (for $n\leq q$):
$$H_n^{\g}(E_{\fin}(\g); \kz) \cong K_{n}(\mZ \g).$$
\end{Thm}

\begin{proof}
We start by recalling
that Farrell-Jones \cite{FJ95} have shown that the Nil-groups
(Bass, Farrell, and Waldhausen type) always
vanish for $n\leq -1$.  So we focus on $n=0,1$.

The proof follows immediately from the discussion above: if one has an infinite
virtually cyclic
subgroup of the form $H\rtimes_\alpha \mathbb Z$ for which the
Farrell Nil-group $NK_n(\mathbb Z H, \alpha)\neq 0$, then the
criterion above tells us that the relative assembly map fails to be
an isomorphism for $n$.  Let us now argue for the converse. We
know that every infinite virtually cyclic subgroup $Q$ of $\g$ is
either a semidirect product $H\rtimes _\alpha \mZ$, or surjects onto
the infinite dihedral group $D_\infty$.  For groups of the first
type, we have that the Farrell Nil-group vanish (by assumption). For
groups of the second type, we know that there is an index two
subgroup which is a semi-direct product of a finite group $H$ with
$\mZ$ (where $\mZ$ acts on $H$ via some automorphism $\alpha$). Our
main theorem says that if the Farrell Nil-group $NK_n(\mathbb Z H,
\alpha)$ associated to the semi-direct product is zero (which holds
by hypothesis), then the Waldhausen Nil-group associated to $V$ are
likewise automatically zero.  In particular, if the Farrell
Nil-groups vanish {\it for every} infinite virtually cyclic subgroup
of the form $H\rtimes _\alpha \mathbb Z$, then the Nil-groups vanish {\it for
every} infinite virtually cyclic subgroup $Q$ in $\g$.  The
criterion discussed above now implies that the relative assembly map
is actually an isomorphism.
\end{proof}

\vskip 10pt

\centerline{\bf Acknowledgements.}

\vskip 10pt

The results in this paper expand on a method used in a previous
paper by the authors \cite{LO}.  The technique in that previous
paper (and hence ultimately the technique in the present paper) was
originally suggested to us by Tom Farrell.  We would like to thank
Tom for his suggestion, as well as for his friendly advice throughout
the years.

In addition, we would like to thank J. Grunewald, D. Juan-Pineda, I. Leary,
S. Prassidis, and M. Varisco for helpful comments on a preliminary draft of
this paper.

The work in this paper was partly supported by the National Science Foundation
under grant DMS-0606002.

\vskip 10pt

\section{Proof of the Main Theorem.}

In this section, we will provide a proof of the main theorem. The
proof will make use of two maps to compare the $K$-theory of
$\mathbb Z \g$ with the $K$-theory of $\mathbb Z \Pi$:
\begin{itemize}
\item the maps $\sigma_*:K_i(\mZ\Pi)\rightarrow K_i(\mZ\g)$,
functorially induced by the inclusion $\Pi\hookrightarrow \g$, and
the transfer maps $\sigma^*:K_i(\mZ\g)\rightarrow K_i(\mZ\Pi)$ (see
Farrell-Hsiang \cite{FH78}). This will be used to establish
$(A)\Rightarrow (B)$.
\item a map $\pi_i\mathcal A: \pi_i \mathcal P(E ; \rho_E)
\longrightarrow  \pi_i \mathcal P(E)$ between suitably defined
spectra of stable pseudo-isotopies (see Farrell-Jones
\cite{FJ95}). This will be used to establish $(B)\Rightarrow (A)$.
\end{itemize}
Recall that in the situation we are dealing with, the group $\Pi =
H\rtimes _\alpha \mathbb Z$ is the canonical index two subgroup of
the group $\g = G_1*_H G_2$.

Another result we will need is that, as shown by Farrell-Hsiang
\cite{FH68}, the group $K_i(\mZ\Pi)$ can be expressed in the
following form:
\begin{equation}
K_i(\mathbb ZH_{\alpha}[\mathbb Z]) \cong C \oplus NK_i(\mathbb ZH,
\alpha) \oplus NK_i(\mathbb ZH, \alpha^{-1})
\end{equation}
where $C$ is a suitable quotient (determined by the automorphism
$\alpha$) of the group $K_{i-1} (\mathbb ZH) \oplus K_i(\mathbb ZH)$.
On the other hand, a classic result of Waldhausen \cite{W1},
\cite{W2} (as reformulated by Connolly-Prassidis \cite{CP02})
expresses $K_i(\mZ\g)$ as:
\begin{equation}
K_i(\mZ[G_1 \ast_{H} G_2]) \cong NK_i(\mathbb ZH;\mathbb Z[G_1 - H],
\mathbb Z[G_2 - H])
\end{equation}
$$ \hskip 1.5in \oplus \big[ \big(K_i(\mZ
G_1)\oplus K_i(\mZ G_2)\big)/K_i(\mZ H)\big].$$

We will first establish that $(A)\Rightarrow (B)$ (in Section 2.1).
We will then briefly recall a construction of Farrell-Jones
\cite{FJ95} of a stratified fiber bundle $\rho_E:E\rightarrow X$,
for a suitably defined space $E$, and stratified control space $X$
(Section 2.2).  We will also explain in that section the relevance
of their result to what we are trying to establish.  Next we shall
use the topology of the spaces $\rho_E^{-1}(x)$ for $x\in X$,
to gain information on the $E^2$-terms in the Quinn
spectral sequence (see Section 2.3).  Finally, we shall use the
information we obtain concerning the spectral sequence to prove that
$(B)\Rightarrow (A)$ (Section 2.4).  Throughout Sections 2.1-2.4, we will
be working exclusively with the case $i=1$.   We will complete the
proof in Section 2.5 by explaining the required modifications needed
to obtain the case $i=0$.

\subsection{Vanishing of Waldhausen Nil forces vanishing of Farrell
Nil}

In order to show that $(A)\Rightarrow (B)$, we first assume that
the Waldhausen Nil-group $NK_1(\mathbb ZH;\mathbb Z[G_1 - H],
\mathbb Z[G_2 - H])=0$.  Note that under this hypothesis, the
decomposition in equation (2) reduces to:
$$K_1(\mZ[G_1 \ast_{H} G_2])\cong \big(K_1(\mZ
G_1)\oplus K_1(\mZ G_2)\big)/K_1(\mZ H)$$ We now want to argue that
the corresponding Farrell Nil-group $NK_1(\mathbb ZH, \alpha)$ is
trivial.  The key observation is contained in the following:

\begin{lemma}
The map $\sigma_*$ is {\it injective} on the subgroup $NK_1(\mathbb
ZH, \alpha)$ in the decomposition (1) of the group $K_1(\mZ\Pi)$.
\end{lemma}

\begin{proof}
This follows from \cite[Proposition 20]{FH68}, which shows that  the
composite map $\sigma ^*\circ \sigma_*:K_1(\mZ\Pi)\rightarrow
K_1(\mZ\Pi)$ has an explicit expression in terms of the
decomposition of $K_1(\mZ\Pi)$ given in equation (1) above: it maps
an element in the $NK_1(\mathbb ZH, \alpha)$ summand to the sum of
itself with the image of this element under the canonical
isomorphism $\tau: NK_1(\mathbb ZH, \alpha)\rightarrow NK_1(\mathbb
ZH, \alpha^{-1})$.  This can be seen as follows: in the short exact
sequence
$$0\rightarrow \Pi \rightarrow \g \rightarrow \mZ /2 \rightarrow 0$$
the $\mZ/2$ acts on $\Pi = H \rtimes _\alpha \mZ$ via the map
$z\mapsto -z$ on the $\mZ$ factor.  But the Farrell Nil-groups
$NK_1(\mathbb ZH, \alpha)$ and $NK_1(\mathbb ZH, \alpha^{-1})$ are
canonically associated to the sub-semirings $\mZ [H\rtimes _\alpha \mZ^+]$
and $\mZ [H\rtimes _\alpha \mZ^-]$, where $\mZ^+, \mZ^-$ refers to the
non-negative and non-positive integers respectively.  Since the
action of $\mZ /2$ on $\Pi$ interchanges these two sub-semirings
inside $\mZ \Pi$, the induced action of the non-trivial element
$g\in \mZ/2$ on the $K_1(\mZ \Pi)$ interchanges the two summands
$NK_1(\mathbb ZH, \alpha)$ and $NK_1(\mathbb ZH, \alpha^{-1})$ via
the canonical isomorphism $\tau$ (this map was precisely the one
used by Farrell-Hsiang to see that these two Nil-groups are
isomorphic). Now the composite map $\sigma ^*\circ
\sigma_*:K_1(\mZ\Pi)\rightarrow K_1(\mZ\Pi)$ is given by the
following formula
$$x\mapsto \Sigma _{g\in \mZ/2} \hskip 3pt g_*(x)$$
where $g_*:K_1(\mZ\Pi)\rightarrow K_1(\mZ\Pi)$ is the map induced on
the $K$-theory of $\mZ\Pi$ via the action of $g$ on $\Pi$ (recall
that $\Pi$ is normal in $\g$).  In the situation we are interested
in, the discussion above implies that for $x\in NK_1(\mathbb ZH,
\alpha)$, we have $(\sigma ^*\circ \sigma_*)(x) = x + \tau(x)$,
where $\tau:NK_1(\mathbb ZH, \alpha) \rightarrow NK_1(\mathbb ZH,
\alpha^{-1})$ is the canonical isomorphism.  This implies that the
composite map $\sigma ^*\circ \sigma_*$ is injective on the
$NK_1(\mathbb ZH, \alpha)$ summand, and hence the map $\sigma_*$
must likewise be injective, concluding the proof of the lemma.
\end{proof}

In particular, in the decomposition of the group $K_1(\mZ[G_1
\ast_{H} G_2])$, we have that the group $NK_1(\mathbb ZH, \alpha)$
injects into the group $\big(K_1(\mZ G_1)\oplus K_1(\mZ
G_2)\big)/K_1(\mZ H)$. Since the latter group is a finitely
generated abelian group, this implies that the group $NK_1(\mathbb ZH,
\alpha)$ is finitely generated. But Grunewald \cite[Theorem 2.5]{G1}
and Ramos \cite{R} independently showed that the groups
$NK_1(\mathbb ZH, \alpha)$ are either trivial or infinitely
generated.  This forces $NK_1(\mathbb ZH, \alpha)$ to vanish,
establishing $(A)\Rightarrow (B)$.

\subsection{Some results of Farrell-Jones}

Farrell and Jones in \cite[Section 2]{FJ95} defined a $3$-dimensional
stratified control space $X$, and constructed from the pair
$\Pi\hookrightarrow \g$ a stratified fiber bundle $E$ over $X$. The
stratified control space $X$ contains a distinguished point $p$, and
the stratified fiber bundle $\rho_E:E\rightarrow X$ in their
construction has the following properties:
\begin{itemize}
\item E is a closed manifold with $\pi_1(E)\cong \g$,
\item for every $x\in X$ satisfying $x\neq p$, the group $\pi_1(E_x)$
is a finite group,
\item for the distinguished point $p\in X$, $\pi_1(E_p)\cong \Pi$.
\end{itemize}
where $E_x=\rho_E^{-1}(x)$.

The important result for our purposes is \cite[Theorem 2.6]{FJ95},
establishing that the group homomorphism
\begin{equation}
\pi_{i}\mathcal A: \pi_i \mathcal P(E ; \rho_E)
\longrightarrow  \pi_i \mathcal P(E)
\end{equation}
is an epimorphism for every integer $i$. Here $\mathcal P(E)$
is the spectrum of stable topological pseudoisotopies on $E$,
$\mathcal P(E, \rho_E)$ is the spectrum of those stable
pseudoisotopies which are controlled over $X$ via $\rho_E$, and
$\mathcal A$ is the `assembly' map.

Now the relevance to the situation we are considering is that, by
results of Anderson and Hsiang \cite[Theorem 3]{AH77}, the lower homotopy
groups of the pseudo-isotopy spectrum coincide (with a shift in
dimension) with the lower algebraic $K$-theory of the integral group
ring of the fundamental group of the space.  In particular, we have
that the right hand side of the map in (3) is given by:
\[
\pi_j \mathcal P(E)=
\begin{cases}
Wh( \g) , & j=-1 \\
\tK(\mathbb Z\g), & j=-2 \\
K_{j+2}(\mathbb Z\g), & j \leq -3
\end{cases}
\]
To understand the left hand side of the map given in (3), we recall that
Quinn \cite[Theorem 8.7]{Qu82} constructed a spectral sequence
$E^n_{s,t}$ abutting to $ \pi_{s+t} \mathcal P(E ; \rho_E)$
with $E^2_{s,t} = H_s(X ; \pi_t \mathcal P(\rho_E))$. Here
$\pi_q \mathcal P(\rho_E)$, $q \in \mathbb Z$, denotes the
stratified system of abelian groups over $X$ where the group above
$x \in X$ is $\pi_q \mathcal P(\rho^{-1}_E(x))$.

Note that by Anderson and Hsiang's result (see  \cite[Theorem
3]{AH77}) mentioned above, we also have that for every $x$:
\[
\pi_j \mathcal P(E_x)=
\begin{cases}
Wh( \pi_1(E_x)) , & j=-1 \\
\tK(\mathbb Z\pi_1(E_x)), & j=-2 \\
K_{j+2}(\mathbb Z\pi_1(E_x)), & j \leq -3
\end{cases}
\]
where $E_x=\rho^{-1}_E(x)$.

\subsection{$E^2$-terms in the Quinn spectral sequence}

Let us now assume that $(B)$ holds, and let us analyze the
$E^2$-terms in the Quinn spectral sequence. Recall that in our
situation, we have that the groups $\pi_1(E_x)$ are all finite,
except at one distinguished point $p$ where $\pi_1(E_p)\cong \Pi$.
In particular, since the lower algebraic $K$-groups of the integral group
ring of a finite group are finitely generated, this implies that the
groups $\pi_j \mathcal P(E_x)$ are finitely generated groups,
except possibly over the distinguished point $p$.

We now focus on the distinguished point $p$, and recall that we are
assuming that the Farrell Nil-group $NK_1(\mathbb ZH, \alpha)$ is
trivial. Since the Nil-groups $NK_1(\mathbb ZH, \alpha)$ and
$NK_1(\mathbb ZH, \alpha^{-1})$ are canonically isomorphic (see
\cite[Proposition 20]{FH68}), we conclude that the Nil-groups vanish.  
This implies, by the formula for $K_1(\mathbb Z\Pi)$ given in equation 
(1) that $K_1(\mathbb Z\Pi)=C$. Recall that $C$ is a suitable quotient 
(determined by the automorphism $\alpha$) of the group $K_{0} (\mathbb ZH) 
\oplus K_1(\mathbb ZH)$. Since $H$ is finite, we have that $K_0(\mathbb ZH)$ and 
$K_1(\mathbb ZH)$ are finitely generated, and it follows that 
$K_1(\mathbb Z\Pi)=C$ is finitely generated. Since $Wh(\Pi)$ is a quotient of 
$K_1(\mathbb Z\Pi)$, this implies that the Whitehead
group $Wh(\Pi)$ is likewise finitely generated.

From the result of Anderson and Hsiang \cite{AH77}, this tells us
that at the distinguished point $p$, $$\pi_{-1} \mathcal P(E_p)= Wh
(\pi_1(E_p)) = Wh( \Pi)$$ is finitely generated.
Furthermore, by results of Farrell-Jones \cite[Corollary 1.3]{FJ95}
we know that for $j\leq -4$, the groups $$\pi_j \mathcal
P(E_p)=K_{j+2}(\mathbb Z\pi_1(E_p))=K_{j+2}(\mathbb Z\Pi)$$ vanish,
and that the group $\pi_{-3} \mathcal P(E_p)=K_{-1}(\mathbb
Z\pi_1(E_p))=K_{-1}(\mathbb Z\Pi)$ is finitely generated.  We
summarize the information we have so far concerning the homotopy
groups $\pi_j \mathcal P(E_x)$ in the following:

\vskip 5pt

\noindent {\bf Fact 1:}  For the distinguished point $p\in X$, we
have that for $j\leq -1$, the groups $\pi_j \mathcal P(E_p)$ are
finitely generated, except possibly for the case where $j=-2$.  For
all other points $x\neq p$ in $X$, the groups $\pi_j \mathcal
P(E_x)$ are finitely generated for all $j \leq -1$.

\vskip 5pt

Now in terms of the stratified space $X$, let us recall how the
$E^2_{p,q}$ term of the spectral sequence can be computed.  The
space $X$ is a stratified 3-dimensional CW-complex, with six cells.
Furthermore, for points $x,y\in X$ lying in the interior of the same
strata $\sigma^p$, we have that $E_x$ is homeomorphic to $E_y$ (we
will denote this common space by $E_{{\sigma}^{p}}$). In particular,
the stratified system of abelian groups is {\it constant} on the
interior of each cell. Then the $E^2_{p,q}$ term is given by the
homology of the chain complex:
\[
\cdots \rightarrow \bigoplus_{{\sigma}^{p+1}} \pi_q \mathcal
P(E_{{\sigma}^{p+1}}) \rightarrow \bigoplus_{{\sigma}^p} \pi_q
\mathcal P(E_{{\sigma}^{p}}) \rightarrow
\bigoplus_{{\sigma}^{p-1}} \pi_q \mathcal P(E_{{\sigma}^{p-1}})
\rightarrow \cdots
\]
where each sum is over all appropriate dimensional strata in the
decomposition of $X$.  Note that from this chain complex, and the
fact that $X$ is 3-dimensional, we immediately see that
$E^2_{p,q}=0$ as soon as $p\leq -1$ or $4\leq p$.  Similarly, since
all the groups in the chain complex are trivial for $q\leq -4$, we
see that for such values of $q$, $E^2_{p,q}=0$. This gives us:

\vskip 5pt

\noindent {\bf Fact 2:} The only $E^2_{p,q}$-terms that can be
non-zero are those for which $0\leq p\leq 3$ and $-3\leq q$. In
particular, the spectral sequence  collapses (by the $E^4$-stage).

\vskip 5pt

Next we observe that the distinguished point $p\in X$ is actually a
0-cell in the CW-complex structure on $X$.  From Fact 1 above, along
with the chain complex describing the $E^2$-terms, we immediately
obtain the following:

\vskip 5pt

\noindent {\bf Fact 3:}  Within the range $q\leq -1$, all the
$E^2$-terms in the spectral sequence are {\it finitely generated},
with the possible exception of the $E^2_{0,-2}$ term.

\vskip 5pt

\subsection{Vanishing of Farrell Nil forces vanishing of Waldhausen Nil}

Now that we have gathered together information on the $E^2$-terms of
the spectral sequence (assuming statement (B) holds), let us show
that statement (A) must also hold.

Observe that, in an arbitrary spectral sequence, if a term
$E^2_{p,q}$ is finitely generated, then $E^k_{p,q}$ is finitely
generated {\it for all $k\geq 2$}. Furthermore, we have that if a
term $E^2_{p,q}=0$, then $E^k_{p,q}=0$ for all $k\geq 2$. From Fact
2, our sequence collapses by the $E^4$-term, and we see that the
only possible non-zero terms satisfying $p+q=-1$ are $E^4_{0,-1}$,
$E^4_{1,-2}$, and $E^4_{2,-3}$.  Combining Fact 3 with the
observation above, we see that $\pi_{-1} \mathcal P(E ; \rho_E)
\cong E^4_{0,-1} \oplus E^4_{1,-2} \oplus E^4_{2,-3}$ is finitely
generated, and by Farrell-Jones \cite[Theorem 2.6]{FJ95}, this group
{\it surjects} onto $\pi_{-1} \mathcal P(E) = Wh(\g)$.  In
particular, we see that the group $Wh(\g)$ must be finitely
generated, and hence the group $K_1(\mathbb Z\g)$ is finitely
generated. From the decomposition in equation (2), we see that the
Waldhausen Nil-group $NK_1(\mathbb ZH;\mathbb Z[G_1 - H], \mathbb
Z[G_2 - H])$ is a direct summand in $K_1(\mathbb Z\g)$, and hence,
must also be finitely generated.  But Grunewald \cite{G1} has shown
that the Waldhausen Nil-group is either trivial or infinitely
generated.  We conclude that the Waldhausen Nil-group is in fact
trivial, completing the proof that $(B)\Rightarrow (A)$.

\subsection{Modifications for the case i=0}

We now proceed to explain the modifications required to deal with
the case $i=0$.  The argument given above would extend verbatim,
provided we had analogues for $i=0$ for the following key results
(known to hold for $i=1$):
\begin{itemize}
\item the result in \cite[Proposition 20]{FH68}, used in the proof
of Lemma 2.1.
\item the results of Ramos \cite{R} and Grunewald \cite{G1}
establishing that the Waldhausen and Farrell Nil-groups are either
trivial or infinitely generated.
\end{itemize}
To explain why these results extend to $i=0$, we recall some basic
facts concerning the suspension functor from rings to rings.  The {\it cone ring} 
$\Lambda \mathbb Z$ of $\mathbb Z$ is the ring of matrices over $\mathbb Z$ 
such that every column and every row contains only finitely many non-zero 
entries. The suspension ring $\Sigma\z$ is the quotient of $\Lambda\z$ by the 
ideal of finite matrices. For an arbitrary ring $R$ we define the suspension ring 
$\Sigma R=\Sigma\z \otimes R$. The suspension fuctor $\Sigma\z \otimes -$ 
has the key property that $K_i(R) \cong K_{i+1}(\Sigma R)$ for all $i \in \mathbb Z$.  

\vskip 5pt

\noindent {\bf Fact 4:} For $\Pi\leq \Gamma$ as in our Main Theorem, with 
respect to the decomposition 
$$K_0(\mathbb ZH_{\alpha}[\mathbb Z]) \cong C \oplus NK_0(\mathbb ZH,
\alpha) \oplus NK_0(\mathbb ZH, \alpha^{-1})$$
we have for all $x\in NK_0(\mathbb ZH,
\alpha)$ the equality $(\sigma ^*\circ \sigma_*)(x) = x + \tau(x)$,
where $\tau:NK_0(\mathbb ZH, \alpha) \rightarrow NK_0(\mathbb ZH,
\alpha^{-1})$ is the canonical isomorphism.

\vskip 5pt

To establish this fact, we merely note that \cite[Proposition 20]{FH68} holds for 
{\it any} ring $R$, and establishes that, for an abstract ring automorphism 
$\alpha$, we have (1) the $K_1$ of the twisted ring $R_\alpha[T]$ can be computed via
$$K_1(R_\alpha [T])\cong C \oplus NK_1(R, \alpha) \oplus NK_1(R, \alpha ^{-1}),$$
where $C$ is some suitable quotient of $K_0(R) \oplus K_1(R)$,  and 
(2) the groups $NK_1(R, \alpha)$ and $NK_1(R, \alpha ^{-1})$ are canonically 
associated to the sub-semirings $R_\alpha[T^+]$ and $R_\alpha [T^-]$.  Here $T$
is an infinite cyclic group, $T^+$ and $T^-$ the semigroup generated by the generator
and the inverse of the generator respectively.

Now observe that there is a canonical ring 
isomorphism $\Sigma (R_\alpha[T]) \cong (\Sigma R)_{1 \otimes \alpha}[T]$ (see \cite[Remark 2.12 part (2)]{G2}).
This induces a commutative diagram:
$$
\xymatrix{K_0(R_\alpha[T]) \ar[rr] \ar[d]_{\cong} & & K_0(R) \ar[d]^{\cong} \\ 
K_1(\Sigma (R_\alpha [T]))  \ar[rr] & & K_1(\Sigma R)
}
$$
which implies that the kernel of the top row is isomorphic to the kernel of the 
bottom row.  But these groups are, by definition, $NK_0(R,\alpha)\oplus NK_0(R,\alpha ^{-1})$ 
and $NK_1(\Sigma R, 1 \otimes \alpha)\oplus NK_1(\Sigma R, 1 \otimes \alpha^{-1})$
respectively.  Furthermore, the functor $\Sigma$ maps the sub-semiring 
$R_\alpha[T^+]$ (and $R_\alpha [T^{-}]$) of the ring $R_\alpha[T]$ to the sub-semiring
$(\Sigma R)_{1 \otimes \alpha}[T^+]$ (and $(\Sigma R)_{1 \otimes \alpha}[T^-]$, respectively)
of the ring $(\Sigma R)_{1 \otimes \alpha}[T]$.  This forces the isomorphism above
to restrict to isomorphisms $NK_0(R,\alpha)\cong NK_1(\Sigma R, 1 \otimes \alpha)$ and
$NK_0(R,\alpha ^{-1})\cong NK_1(\Sigma R, 1 \otimes  \alpha^{-1})$ respectively.
Finally, under the situation we are dealing with, we note that the non-trivial element
$g\in \mathbb Z/2 = \Gamma/\Pi$ which acts on $R = \mathbb Z\Pi$ by interchanging
the two sub-semirings $R_\alpha[T^+]$ and $R_\alpha [T^-]$ commutes with the functor
$\Sigma$, and hence acts on $(\Sigma R)_{1 \otimes \alpha}[T]$ by permuting the
two sub-semirings $(\Sigma R)_{1 \otimes \alpha}[T^+]$ and $(\Sigma R)_{1 \otimes \alpha}[T^-]$.
This implies that the map $\tau$ acts on the decomposition of $K_0(\mathbb Z\Pi)$ 
by interchanging the two factors $NK_0(R,\alpha)$ and $NK_0(R,\alpha ^{-1})$, as 
was required in the argument of Lemma 2.1.  This completes the verification of the 
first point mentioned above.

For the second point mentioned above, we note that, by the argument in the previous 
paragraph, we have that $NK_0(\mathbb ZH, \alpha)  \cong NK_1(\Sigma (\mathbb ZH), 1 \otimes
\alpha)$.  Since Grunewald \cite{G1} showed that the $NK_1(R, \beta)$ is either 
trivial or infinitely generated, for {\it any} ring $R$ and {\it any} automorphism $\beta$
of finite order, 
we immediately have the corresponding property for $NK_0$.  A similar argument 
shows that, for the Waldhausen Nil-groups, we have isomorphisms:
$$NK_0(\mathbb ZH;\mathbb Z[G_1 - H],\mathbb Z[G_2 - H]) 
\cong NK_1\big(\Sigma(\mathbb ZH); \Sigma(\mathbb Z[G_1 - H]),\Sigma(\mathbb Z[G_2 - H])\big)$$
Note in the above expression that the functor $\Sigma$ has a natural extension to 
left bimodules, in the sense that if $B$ is a left bimodule for the ring $R$, then $\Sigma B$
is a left bimodule for $\Sigma R$ (see \cite[Section 2]{G2} for more details).  Applying 
Grunewald's result for $NK_1$ to the right hand term, 
concludes the argument for the following:

\vskip 5pt

\noindent {\bf Fact 5:} For $\Pi\leq \Gamma$ as in our main theorem, we have that both
\begin{itemize}
\item the Farrell Nil-group $NK_0(\mathbb ZH, \alpha)$, and
\item the Waldhausen
Nil-group $NK_0(\mathbb ZH;\mathbb Z[G_1 - H],\mathbb Z[G_2 - H])$,
\end{itemize}
are either trivial or infinitely generated.

\vskip 5pt

Finally, we explain how Facts 4 and 5 can be used to promote the arguments in 
Sections 2.1-2.4.  Using Fact 4, the proof in Lemma 2.1 extends verbatim to show
that the map $\sigma_*: K_0(\mathbb Z\Pi)\rightarrow K_0(\mathbb Z\Gamma)$ is 
injective on the subgroup $NK_0(\mathbb ZH, \alpha)$.  Using Fact 5 and the 
decompositions of $K_0(\mathbb Z\Pi)$ and $K_0(\mathbb Z\Gamma)$ given in
equations (1) and (2), the arguments in the last paragraph of Section 2.1 yield
the implication $(A) \Rightarrow (B)$ for the case $i=0$.

In the proof of the converse, the sole changes occur in Section 2.3, Fact 1, 
where all the groups
$\pi_j \mathcal P(E_p)$ are finitely generated except possibly for $j=-1$ (instead of
$j=-2$).  Correspondingly, there is a change in Fact 3, where within the range $q\leq -1$,
all the $E^2$-terms of the spectral sequence are finitely generated except possibly for
the $E^2_{0,-1}$ term.  The argument in Section 2.4 allows one to conclude that the
group $\pi_{-2} \mathcal P(E ; \rho_E) \cong E^4_{0,-2} \oplus E^4_{1,-3}$ is finitely
generated, and the rest of the discussion in Section 2.4, along with 
Fact 5, completes the converse implication $(B)\Rightarrow (A)$ for the case $i=0$.

\section{Concluding remarks.}

While our main theorem provides a relationship between the Farrell
Nil-groups and the Waldhausen Nil-groups, one can speculate about
whether one can establish a further reduction to the Bass
Nil-groups.  This motivates the following:

\vskip 5pt

\noindent {\bf Question:} For an arbitrary ring $R$, are the
following two statements equivalent:
\begin{itemize}
\item the Bass Nil-group $NK_i(R)$ vanishes, and
\item the Farrell Nil-groups $NK_i(R, \alpha)$ vanish for {\it
every} automorphism $\alpha$ of the ring $R$.
\end{itemize}

\vskip 5pt

This last question seems extremely difficult to answer.  The most
interesting case would be where the ring $R$ is the integral group
ring of a finitely generated group $G$, and the automorphisms are
generated by automorphisms of $G$.

Finally, we point out that if one does {\it not} assume vanishing of
the Farrell Nil-group of $\Pi$, the Quinn spectral sequence still
gives us some information relating the Farrell Nil-group of $\Pi$
with the Waldhausen Nil-group of $\g$.  Indeed, a more careful
analysis of the Quinn spectral sequence can be used to show that the
Farrell Nil-group appearing in the $E^2$ term survives to the
$E^4$ stage, and {\it surjects} onto the Waldhausen Nil-group (see
Prassidis \cite[pages 412-413]{P97} for a similar analysis).  This is
motivated by the well-known philosophy that, in the map of
Farrell-Jones:
$$
\pi_{i}\mathcal A: \pi_i \mathcal P(E ; \rho_E)
\longrightarrow  \pi_i \mathcal P(E)
$$
the controlled (resp. non-controlled) part of the pseudo-isotopy
spectrum $\mathcal P(E ; \rho_E)$ must map to the controlled
(resp. non-controlled) part of the pseudo-isotopy spectrum $\mathcal
P(E)$.  The non-controlled part is precisely the
Nil-groups.

Now one could use this surjection to show that the exponent of the
Waldhausen Nil-group must divide the exponent of the Farrell
Nil-group.  We observe however that the estimates this would yield
provide no improvement on known estimates (due independently to
Grunewald \cite{G1} and Ramos \cite{R}).  In the interest of clarity
of exposition, the authors have chosen to omit this further
analysis, leaving the details to the interested reader.

\end{document}